\documentclass{amsart}
\usepackage{epsfig}
\usepackage[english]{babel}
\usepackage{amssymb}
\usepackage{amsmath}
\usepackage{graphicx}
\newcommand{\1}{\mathbb{1}}

\newcommand{\Z}{\mathbb{Z}}

\newcommand{\ds}{\displaystyle}

\newtheorem{thank}{\ \ \ Acknowledgment}
\newcounter{tictac}

\def\1{\,\rlap{\mbox{\small\rm 1}}\kern.15em 1}

\def\build#1_#2^#3{\mathrel{\mathop{\kern 0pt#1}\limits_{#2}^{#3}}}
\def\tend#1#2{\build\hbox to 12mm{\rightarrowfill}_{#1\rightarrow #2}^{ }}

\def\converge#1#2#3#4{\build\hbox to
#1mm{\rightarrowfill}_{#2\rightarrow #3}^{\hbox{\scriptsize #4}}}

\theoremstyle{definition}

\newtheorem{thm}{Theorem}[section]
\newtheorem{prop}[thm]{Proposition}
\newtheorem{lemm}[thm]{Lemma}

\newtheorem{rem}[thm]{Remark}
\newtheorem{Cor}[thm]{Corollary}
\newcommand{\beq}{\begin{equation}}
\newcommand{\eeq}{\end{equation}}
\newcommand{\xbm}{(X,{\mathcal B},\mu)}

\begin{document}
\title[Approximately transitive dynamical systems
and simple spectrum]{Approximately transitive dynamical systems
and simple spectrum}
%Countable Abelien group action with simple spectrum and without Approximately transitive }
\author{E. H. El Abdalaoui}
\address{ Department of Mathematics, University
of Rouen, LMRS, UMR 60 85, Avenue de l'Universit\'e, BP.12, 76801
Saint Etienne du Rouvray - France}
\email{elhoucein.elabdalaoui@univ-rouen.fr }

\author {M. Lema\'nczyk $^1$}
\address{Faculty of Mathematics and Computer Science,
Nicolas Copernicus University, Chopin street 12/18, 87-100,
Toru\'n, Poland} \email{mlem@mat.uni.torun.pl} \footnote{Research
supported by the EU Program Transfer of Knowledge ``Operator
Theory Methods for Differential Equations" TODEQ and the  Polish
MNiSzW grant N N201 384834}

\maketitle

{\renewcommand\abstractname{Abstract}
\begin{abstract} For some countable discrete torsion Abelian
groups we give examples of their finite measure-preserving actions
which have simple spectrum and no approximate transitivity
property.
\vspace{8cm}\\

\hspace{-0.7cm}{\em AMS Subject Classifications} (2000): 37A25
\end{abstract}
%-------------------------------------------------- Introduction ----------------------------------------------------------
\thispagestyle{empty}
\newpage
\section{Introduction}

The property of approximate transitivity (AT) has been introduced
to ergodic theory by Connes and Woods \cite{Connes-woods} in 1985
in connection with their measure-theoretic characterization of
ITPFI hyperfinite factors in the theory of von Neumann algebras.
Since then the approximately transitive dynamical systems have
been objects of study in several papers, e.g.\
\cite{elabdal-Mariusz}, \cite{Dooley-Quas}, \cite{Gi-Ha},
\cite{Golodets}, \cite{Hawkins}, \cite{HawkinsR}, \cite{Sokhet}.
All systems with the AT property have zero entropy
\cite{Connes-woods}, \cite{Da} and it is only recently
\cite{elabdal-Mariusz}, \cite{Dooley-Quas}  that explicit examples
of zero entropy systems without AT property have been discovered.

In view of the form of its definition (see below) the AT property
seems to be related to spectral properties, or more precisely to
the spectral multiplicity of the measure-preserving system.
Indeed, it was already in the 1990.\ when J.-P.\ Thouvenot using
some generic type arguments observed that the AT property implies
the existence of a cyclic vector in the $L^1$-space of the
underlying dynamical system. Moreover, a certain modification of
the definition of the AT property (in which we replace $L^1$-norms
by $L^2$-norms) which was done in \cite{HawkinsR} implies
$L^2$-simplicity of the spectrum. However, in general it is still
an interesting open problem whether the original AT property
implies simplicity of $L^2$-spectrum, i.e.\ simplicity of the
spectrum of the classical Koopman representation of the dynamical
system. Moreover, as noticed for example in \cite{Dooley-Quas},
neither the other implication was known to hold. The aim of this
note is to answer this latter question -- we give examples of
systems with are not AT but have simple $L^2$-spectrum.

One more motivation for this note is the problem of a
measure-theoretic characterization of simplicity of
$L^2$-spectrum. This problem has a long history in ergodic theory,
see for example \cite{Ka-Th}, \cite{Le} and the references therein
and our examples  bring the answer to one more open question
(stated explicitly in case of $\Z$-actions in \cite{Fe1} as well
as in \cite{Dooley-Quas}): namely we give the negative answer to
the question whether the class of funny rank one systems
\cite{Fe}, \cite{delJunco} coincides with the class of systems
with simple $L^2$-spectrum; indeed a funny rank one system enjoys
the AT property (see below).

We would like however to emphasize that in the note we do not
consider $\Z$-actions which are the most popular objects of study
in ergodic theory, and therefore the problems we mentioned above
stay open for the class of $\Z$-actions (and actions of many other
groups). The systems which are considered below are actions of
some torsion Abelian countable discrete groups. More precisely
$G=\bigoplus_{n=0}^\infty \Z/p_n\Z$ where $(p_n)_{n\geq0}$ is an
increasing sequence of primes numbers. We will then deal with so
called Morse extensions given by some special $\Z/2\Z$-valued
cocycles over a discrete spectrum action of $G$ studied by M.\
Guenais in \cite{Gu} (in that paper she proved that the resulting
$G$-actions have Haar component of multiplicity one in the
$L^2$-spectrum). All such systems have simple spectrum. We will
show that  for each sufficiently fast choice of parameters in
Guenais' construction we can apply the criterion of being non-AT
system formulated in \cite{elabdal-Mariusz} (this criterion is an
elaborated version of the method already used in
\cite{Dooley-Quas}). More precisely, we will prove the following.

\begin{thm}\label{SnoAT}\em Assume that
$G=\bigoplus_{n=0}^\infty \Z/p_n\Z$ with $(p_n)_{n\geq0}$ an
increasing sequence of primes numbers and $p_n\geq 5^{2(n+1)}$,
$n\geq 0$. Then there exists a simple $L^2$-spectrum  action of
$G$ without the AT property.
\end{thm}

\vspace{1ex}

Following Connes and Woods \cite{Connes-woods} we now recall the
definition of the AT property.

Let $G$ be an Abelian countable discrete group. Assume that this
group acts as measure preserving maps: $g\mapsto T_g\in
Aut(X,\mathcal{ B},\mu)$, where $(X,\mathcal{B},\mu)$ is a
standard probability Borel space.  The action
$(X,\mathcal{B},\mu,T)$ with $T=(T_g)_{g\in G}$ (or simply $T$) is
called AT if for an arbitrary family of nonnegative functions
$f_1,\ldots,f_l \in L_{+}^1(X,\mathcal{B},\mu), l \geq 2,$ and any
$\varepsilon
>0$, there exist a positive integer $s$, elements
$g_1,\ldots,g_s \in G$, real numbers $\lambda_{j,k} >0,
j=1,\ldots, l, k=1, \ldots, s$ and $f \in
L_{+}^1(X,\mathcal{B},\mu)$ such that \beq\label{f0}
 \|f_j-
\sum_{k=1}^{s}\lambda_{j,k}f\circ T_{g_k}\|_1 < \varepsilon,~1
\leq j \leq l. \eeq Recall now that $T$ has the {\em funny rank
one} property \cite{Fe}, \cite{delJunco} if for every
$A\in\mathcal{B}$ and $\varepsilon>0$ we can find
$F\in\mathcal{B}$, $g_1,\ldots,g_N\in G$ such that the family
$\mathcal{R}=\{T_{g_1}F,\ldots, T_{g_N}F\}$, called a {\em funny
Rokhlin tower}, consists of sets which are disjoint and for some
$J\subset \{1,\ldots,N\}$
\beq\label{f1}\mu\left(A\triangle\bigcup_{i\in
J}T_{g_i}F\right)<\varepsilon.\eeq Alternatively, $T$ is of funny
rank one if and only if there exists a sequence
$(\mathcal{R}_n)_{n\geq1}$ of funny towers with bases $F_n$ such
that each set $A\in\mathcal{B}$ can be $\varepsilon$-approximated
(as in~(\ref{f1})) by the union of some levels for each
$\mathcal{R}_n$ whenever $n\geq n_\varepsilon$. From this it
easily follows that each system which is of funny rank one is AT.
Indeed, the set of functions which are constant on levels of
$\mathcal{R}_n$  for some $n\geq1$ (and zero outside
$\bigcup\mathcal{R}_n$) are dense in $L^1$, and the set of those
which are additionally positive is dense in $L^1_+$; it follows
that given $f_1,\ldots,f_l\in L^1_+(X,\mu)$ and $\varepsilon>0$ it
is enough to take $f=\1_{F_n}$ (in~(\ref{f0})) for $n\geq1$ large
enough.

The action $T$ of $G$ on $\xbm$ induces a (continuous) unitary
representation, called a {\em Koopman representation}, of $G$ in
the space $L^2(X,\mathcal{B},\mu)$ given by $U_{T_g}f = f\circ
T_g, f \in L^2(X,\mathcal{B},\mu)$ and $g \in G$. We recall that a
Koopman representation is said to have {\em simple spectrum} if
for some $f\in L^2\xbm$, $L^2\xbm=\overline{\mbox{span}}\{f\circ
T_g:\:g\in G\}$. Given $f\in L^2(X,\mu)$ we define its {\em
spectral measure} $\sigma_f$ (or in a more precise notation,
$\sigma_{U_T,f}$) to be a finite Borel measure on the dual
$\widehat{G}$ of $G$ determined by
$$
\widehat{\sigma}_{U_T,f}(g)=
\widehat{\sigma}_f(g):=\int_{\widehat{G}}\chi(g)
\,d\sigma_f(\chi)=\int_{X}f\circ
T_g\cdot\overline{f}\,d\mu,\;\;\mbox{for all}\;g\in G.$$  See
\cite{Ka-Th}, \cite{Le} for more information on the spectral
theory of $G$-actions.

\section{A criterion for a system to be non-AT}

Let $G$ be an Abelian discrete countable group which we assume to
act  on a probability standard Borel space $(X,\mathcal{B},\mu)$
as measure-preserving maps: $g\mapsto T_g$. Assume that
$\mathcal{P}=\{P_0,P_1\}$ a partition of $X$ (with
$P_0\in\mathcal{B}$). Through its $\mathcal{P}$-names every point
$x\in X$ can be now coded: $x\mapsto\pi(x)=(x_g)_{g\in G}$ where
\[
x_g=\left\{
     \begin{array}{ll}
       0 & \hbox{if $T_g(x) \in P_0$} \\
       1 & \hbox{if not.}
     \end{array}
   \right.
\]

Let $\Lambda$ be a finite subset of $G$. By a {\em funny word on
the alphabet $\{0,1\}$ based on} $\Lambda$ we mean a sequence
$(W_g)_{g \in \Lambda}$ with $W_g \in \{0,1\}$, $g\in\Lambda$. For
any two funny words $W,W'$ based on the same set $\Lambda\subset
G$ their {\em Hamming distance} is given by
\[
\overline{d}_{\Lambda}(W,W')=\frac1{|\Lambda|}\left|\left \{ g \in
\Lambda:~W_g \neq W'_g \right \}\right|.
\]

As noticed in \cite{elabdal-Mariusz} we then have the following
extension of Dooley-Quas'  necessary condition for a system
$T=(T_g)_{g\in G}$ to have the AT property  \cite{Dooley-Quas}.

\begin{prop}{(\cite{Dooley-Quas})}\label{fanny} \em
Let $(X,\mathcal{B},\mu,T)$ be an AT dynamical system. Then for
any $\varepsilon >0$ there exist a finite set  $\Lambda \subset G$
and a funny word $W$  based on $\Lambda$ such that
\[
|\Lambda|\mu\left(\{x\in X:~
\overline{d}_\Lambda(\pi(x)|_{\Lambda},W)<\varepsilon\}\right)>
1-\varepsilon.
\]
\end{prop}

\vspace{1ex}

The contraposition in Proposition~\ref{fanny} yields automatically
a criterion for a system to be non-AT. Some further work has been
done in \cite{elabdal-Mariusz} to formulate a condition stronger
than the negation of the assertion in Proposition~\ref{fanny} and
which may be applied to many systems (the criterion obtained this
way is an elaborated version of the method already used in
\cite{Dooley-Quas}). We now present this criterion
(Proposition~\ref{noat} below) in its generality needed for this
note.

A probability Borel measure $\rho$ defined on  $\widehat{G}$ is
called a {\em strong Blum-Hanson measure} (SBH measure shortly) if
the following holds
\[
\limsup_{k \to+\infty}{\sup_{\Theta \subset G, |\Theta|=k
,(\eta_\theta)_{\theta \in \Theta}\in \{1,-1\}^k}}
\int_{\widehat{G}}\left| \frac1{\sqrt{k}} \sum_{\theta \in \Theta}
\eta_{\theta}\cdot\chi(\theta) \right |^2 d\rho(\chi) <2.
\]
\

\noindent{}Clearly the Haar measure $m_{\widehat{G}}$ of
$\widehat{G}$ is an SBH measure; more generally, every absolutely
continuous probability measure with density $d\in
L^1(\widehat{G},m_{\widehat{G}})$ satisfying $\|d\|_{\infty}<2,$
is an SBH measure. It is shown in \cite{elabdal-Mariusz} that each
SBH measure is a Rajchman measure (i.e.\ its Fourier transform
vanishes at infinity).

\begin{prop}[\cite{elabdal-Mariusz}]\label{noat}\em
Assume that a dynamical system $(X,\mathcal{B},\mu,T)$ is ergodic
and that there exists a partition $\mathcal{P}=\{P_0,P_1\}$ with
the following properties:
\begin{enumerate}
    \item[i)] There exists $S\in Aut(X,\mathcal{B},\mu)$
    commuting with all elements $T_g$, $g \in G$, such that
    $SP_0=P_1$.
    \item[ii)] The spectral measure
    $\sigma_{\1_{P_0}-\1_{P_1}}$
    is an SBH measure.
\end{enumerate}
Then the system $T=(T_g)_{g\in G}$  is not AT.
\end{prop}

\section{Countable Abelian discrete group
action with simple spectrum and without the AT property}

In this section we shall present the proof of Theorem~\ref{SnoAT}.

\subsection{Cocycles for discrete spectrum
actions of countable discrete Abelian groups} For completness, we
now briefly present an extension of the classical spectral theory
of compact group extensions  of $\Z$-actions  to group actions.

Let $G$ be a countable Abelian discrete group acting on a standard
probability Borel space $(X,\mathcal{B},\mu)$ as
measure-preserving maps: $g\mapsto T_g\in Aut(X,\mathcal{B},\mu)$.
Assume that $\Gamma$ is a compact metric Abelian (written
multiplicatively) group with Haar measure $m_\Gamma$. A
{\it{cocycle}} associated to $T=(T_g)_{g\in G}$ with values in
$\Gamma$ is a measurable function $ \varphi:X \times G \to\Gamma$
which satisfies
\beq\label{coc}\varphi(x,g+g')=\varphi(x,g)\varphi(T_gx,g')\eeq
 for a.e.\ $x\in X$ and all $g,g'\in G$.  We will also write
$\varphi_g(x)$ instead of $\varphi(x,g)$. Then the $\Gamma$-{\em
extension of $T$ associated to} $\varphi$ is the $G$-action
$T_{\varphi}=\left((T_\varphi)_g\right)_{g\in G}$ defined on $(X
\times \Gamma,\mu \otimes m_\Gamma)$ by
$$
(T_{\varphi})_g:X \times \Gamma  \to X \times \Gamma,\;\;
  (x,\gamma) \mapsto (T_gx,\varphi_g(x)\gamma).
$$
The space $L^2(X \times \Gamma , \mu \times m_\Gamma)$ can be
decomposed as \beq\label{d1} L^2(X \times \Gamma , \mu \times
m_\Gamma )=\bigoplus_{\chi \in \widehat{\Gamma}}L_{\chi}, \eeq
where each of the subspaces $L_{\chi} := \{f \bigotimes \chi : f
\in L^2(X, \mu)\}$ is ${(U_{{T_{\varphi}}})}_g$-invariant, and the
restriction of the Koopman representation $U_{T_{\varphi}}$ of $G$
to $L_{\chi} $ is, via the map \beq\label{d0} f\otimes \chi\mapsto
f,\;f\in L^2(X,\mu),\eeq unitarily equivalent to the
$G$-representation
$V_{\varphi,T,\chi}=\left(V_{\varphi,T,\chi,g}\right)_{g\in G}$
given by
$$
V_{\varphi,T,\chi,g}:L^2(X, \mu)\to L^2( X , \mu ),\;\;
V_{\varphi,T,\chi,g}(f)(x) = \chi(\varphi_g(x)) f (T_gx) , x \in
X.$$ It follows that the spectral properties of $T_{\varphi}$ are
determined by the spectral properties of $V_{\varphi,T,\chi}$,
$\chi \in \widehat{\Gamma}$.

Assume now that $\Gamma=\Z_2=\{1,-1\}$. Set
$$
S:X \times \Z_2\to X\times \Z_2,\;\;
S(x,\varepsilon)=(x,-\varepsilon)$$ and notice that $S$ commutes
with all automorphisms $(T_\varphi)_g$, $g\in G$, i.e.\ $S$ is in
the centralizer of $T_{\varphi}$. In this case the
decomposition~(\ref{d1}) is given by $L^2(X \times \Z_2,\mu\otimes
m_{\Z_2})=L_0 \oplus L_1$ where $L_0=\{F \in L^2(X \times
\Z_2,\mu\otimes m_{\Z_2}):\:F\circ S=F \}$ and $L_1=\{F\in L^2(X
\times \Gamma,\mu\otimes m_{\Z_2}):\:F\circ S=-F\}$. Thus the
restriction of the Koopman representation $U_{T_{\varphi}}$ on
$L_1$ is unitarily equivalent to the representation $V=(V_g)_{g
\in G}$ on $L^2(X,\mu)$ given by
\[
(V_g(f))(x)={\varphi_g(x)} f(T_gx),~{\rm {for ~all~}} x \in X,\; g
\in G.
\]
Let $\mathcal{P}=\{P_0,P_1\}$ be the partition of $X\times\Z_2$
given by $P_0=X \times \{1\}, P_1=X \times \{-1\}$. Then
$\1_{P_0}-\1_{P_1}\in L_1$, in fact
$$
\1_{P_0}-\1_{P_1}=\1_X\otimes\left(\1_{\{1\}}-\1_{\{-1\}}
\right)=\1_X\otimes \chi$$ where $\chi$ is the only non-trivial
character of $\Z_2$; in particular
$\|\1_{P_0}-\1_{P_1}\|_{L^2}=1$, so its spectral measure is a
probability measure. In view of~(\ref{d0}) we have
$$
\widehat{\sigma}_{U_{T_\varphi},\1_{P_0}-\1_{P_1}}(g)=
\widehat{\sigma}_{V_{\varphi,T,\chi},\1_X}(g)=$$$$ \langle
V_{\varphi,T,\chi,g}\1_X,\1_X\rangle_{L^2(X,\mu)}=
\int_X\chi(\varphi_g(x))\,d\mu(x) =\int_X\varphi_g(x)\,d\mu(x).$$
Moreover, $SP_0=P_1$. In this way we have proved the following.

\begin{lemm}\label{pomoc} \em For any ergodic $G$-action
$T=(T_g)_{g\in G}$  and its $\Z_2$-extension $T_\varphi$ for the
partition $\mathcal{P}=\{P_0,P_1\}$ defined above we have:

(i) The element $S$ defined above is in $Aut(X\times\Z_2,
\mu\otimes m_{\Z_2})$ and belongs to the centralizer of
$U_{T_\varphi}$.

(ii) The spectral measure $\sigma_{\1_{P_0}-\1_{P_1}}=
\sigma_{U_{T_\varphi},\1_{P_0}-\1_{P_1}}$ satisfies \beq\label{d3}
\widehat{\sigma}_{\1_{P_0}-\1_{P_1}}(g)=
\int_X\varphi_g(x)\,d\mu(x)\;\;\mbox{for all}\;g\in
G.\eeq\end{lemm}

\vspace{1ex}

Comparing the above lemma with Proposition~\ref{noat} we can see
that the only thing which is missing  is the fact that in general
the measure $\sigma_0=\sigma_{U_{T_\varphi},\1_{P_0}-\1_{P_1}}$ is
not SBH, and to achieve non-AT property we will have to control
the Fourier transform of $\sigma_0$ to obtain absolute continuity
and a relevant boundedness of the density.

Recall also that if in addition the spectrum of the $G$-action $T$
is discrete then by an obvious extension of Helson's analysis from
\cite{Helson} we obtain that the spectral type of $V$ has the
purity law: it is either discrete, or continuous and purely
singular, or equivalent to Haar measure $m_{\widehat{G}}$.

\subsection{Morse cocycles and Guenais' constructions}

We now present Guenais' construction from   \cite{Gu}. Assume that
$(p_n)_{n\geq0}$ is an increasing sequence of prime numbers and
let $G=\oplus_{n \geq 0}{\Z}/{p_n\Z}$. Consider the action of $G$
on its dual  $X=\prod_{n=0}^{+\infty} {\Z}/{p_n\Z}$ by
translations $T_gx=x+g$, where the space $X$ is equipped with Haar
measure $\mu=m_{X}$ (which is the product of counting measures
$m_{\Z/p_n\Z}$ on coordinates). The resulting $G$-action has
discrete spectrum. Moreover, the action $T=(T_g)_{g\in G}$ has
funny rank one. Indeed, to see it set
$$F_n=\{x \in X:\: x_0=\cdots=x_{n-1}=0\},\;
G_n=\bigoplus_{k=0}^{n-1}{\Z}/{p_n\Z},\,n\geq0.$$
 The sequence $(\mathcal{R}_n)$ of funny
Rokhlin towers (see~(\ref{f1})) is given by $(T_gF_n)_{g \in G_n}$
where $F_0=X$, $G_0=\{0\}$. Note that the union of levels of each
such Rokhlin tower fills up the whole space.

As before we consider $\Gamma=\Z_2$ and by a {\em Morse cocycle}
one means a cocycle which will be constant on all levels of the
funny Rokhlin towers  $(T_gF_n)_{g \in G_n}$, $n\geq1$, in the
sense made precise below.

Take an arbitrary sequence of $(\varepsilon_n(j))_{0 \leq j \leq
p_n}, n \geq 0)$ with values in $\Z_2$ in which
$\varepsilon_n(0)=1$. Then define
\begin{eqnarray}\label{Mcocycle}
\varphi(x,g)=\varphi_g(x)=\prod_{n=0}^{\infty}\varepsilon_n(x_n)
\varepsilon_n(x_n+g_n) ~~{\rm {for~all~~}} (x,g) \in X \times G.
\end{eqnarray}
By (\ref{Mcocycle}) it follows directly that the cocycle
identity~(\ref{coc}) is satisfied and moreover for $x\in F_n$ and
$g\in G_n$ we have $$\phi_g(x)=\prod_{k=0}^\infty
\varepsilon_k(g_k)=\prod_{k=0}^{n-1}\varepsilon_k(g_k),$$ so
$\phi_g$ is constant on $F_n$ for $g \in G_n$
($\phi_g(x)=\phi_g(0)$) and in fact since
$$ \varphi(T_hx,g)\varphi(x,g)=\varphi(x,g+h),$$ we also have
$\varphi_g$ is constant on each level $T_hF_n$, $h\in G_n$,
whenever $g\in G_n$: $\phi_g|_{T_hF_n}=\phi_{g+h}(0)\phi_g(0)$.
(In fact, as shown in \cite{Gu}, each Morse cocycle is defined by
a $\Z_2$-valued sequence $(\varepsilon_n(j))_{0 \leq j \leq p_n},
n \geq 0)$.) The following has been observed in \cite{Gu}:
\beq\label{mult} \mbox{For each Morse cocycle $\varphi$, the
representation $V$ on $L_1$ has simple spectrum.}\eeq

Now by the purity law it follows that once the spectral measure of
$\1_X$ (for $V$) is continuous, the spectrum on $L_1$ is
continuous and therefore the Koopman representation
$U_{T_\varphi}$ has simple spectrum whenever $\varphi$ is a Morse
cocycle.

Set $ \varepsilon_n(k)=  \ds \left(\frac{k}{p_n}\right)$ for $0 <
k < p_n$ and $\varepsilon_n(0)=1$, where $\ds
\left(\frac{k}{p_n}\right)$ is the Legendre symbol, that is, it
is equal to~1 if $k\neq0$ is a square modulo~$p_n$, $-1$ if not
and $\ds \left(\frac{0}{p_n}\right)=0$.

\begin{thm}[\cite{Gu}] \label{Riesz}\em The spectral
measure $\sigma_0=\sigma_{V,\1_X}$ of the function $\1_X$ is the
product measure $\bigotimes_{n=0}^\infty|P_n|^2m_{\Z/p_n\Z}$ where
$P_n$ is defined on ${\Z}/{p_n\Z}$ by
$$P_n(x)=\frac1{\sqrt{p_n}}\left (1+\sum_{k=1}^{p_n-1}
\left(\frac{k}{p_n}\right) e^{-2i\pi \frac{k x}{p_n}}\right).$$
Moreover, $\sigma_0$ is equivalent to the Haar measure of the dual
group $\widehat{G}$ with \beq\label{mel1}
\left(\bigotimes_{n=0}^\infty|P_n|^2m_{\Z/p_n\Z}\right)
\widehat{}\;(g)=\widehat{\sigma}_0(g)\;\;\mbox{for each}\;g\in
G.\eeq
\end{thm}

\vspace{1ex}

It follows that the sequence
$\left(\bigotimes_{n=0}^N|P_n|^2m_{\Z/p_n\Z}\right)_{N\geq1}$ of
measures converges weakly to $\sigma_0$ and since $\sigma_0\ll
m_X$ we obtain the following.

\begin{Cor}\label{mel2}\em  The sequence
$\left(\bigotimes_{n=0}^N|P_n|^2\right)_{N\geq1}$ converges weakly
in $L^1(X,m_X)$ to $\frac{d\sigma_0}{d m_X}$.\end{Cor}

\vspace{1ex}

Notice that the sequence $(|P_n|^2)_{n\geq0}$ meant as polynomials
on $X$ ($P_n(x)=P_n(x_n)$ for $x\in X$) is {\em ultra flat}, that
is
$$
\frac{\|P_n|\|_{L^2}}{\|P_n\|_{L^\infty}}\to 1;$$ indeed,
$P_n(0)=\frac1{\sqrt{p_n}}$ and for $x\neq0$ by the Gauss formula
\cite{BerndtEW}, \cite{Rosen}
$$P_n(x)=\frac1{\sqrt{p_n}}\left (1+\sum_{k=1}^{p_n-1}
\left(\frac{k}{p_n}\right) e^{-2i\pi \frac {k x}{p_n}}\right)=
$$$$
\frac1{\sqrt{p_n}}\left (1+\sum_{k=0}^{p_n-1}e^{-2i\pi k^2x/p_n} -
\sum_{k=0}^{p_n-1}e^{-2i\pi kx/p_n}\right)=\frac1{\sqrt{p_n}}
\left(1+\delta_{p_n}\sqrt{p_n}\left(\frac{x}{p_n}\right)\right)$$
with $\delta_{p_n}=1$ if $p_n\equiv 1$ mod~$4$ or $\delta_{p_n}=i$
if $p_n\equiv 3$ mod~$4$; whence \beq\label{mel3}
1-\frac1{\sqrt{p_n}}\leq \left|P_n(x)\right|\leq
1+\frac1{\sqrt{p_n}}.\eeq

Recall also the following elementary fact.
\begin{lemm}\label{ho}\em Let $a \in (0,1)$.
Then we have $$\ds \prod_{n=1}^{\infty}(1+a^n) \leq \exp\left({\ds
\frac{a}{1-a}}\right).$$
\end{lemm}

\vspace{2ex}

In view of (\ref{mel3}) and Corollary~\ref{mel2} \beq\label{mel4}
\frac{d\sigma_0}{dm_X}\leq\prod_{n=0}^\infty
\left(1+\frac1{\sqrt{p_n}}\right)^2.\eeq

Let us choose now the sequence $(p_n)_{n\geq0}$ so that
$\sqrt{p_n} \geq 5^{n+1}$ for any $n \geq 0$. By Lemma~\ref{ho} we
have
$$ \ds \prod_{n=0}^{\infty}( 1+\frac1{\sqrt{p_n}}) \leq e^{1/4}.$$
Hence by (\ref{mel4}), $\frac{d\sigma_0}{dm_X}\leq e^{1/2}<2$ and
the proof of Theorem~\ref{SnoAT} is complete.

\begin{rem}The action of the group
$G=\bigoplus_{n=0}^\infty \Z/p_n\Z$ on its dual by translations
(as above) is a particular case of so called $(C,F)$-actions of
the group $G$  (see \cite{Dan1}), all such actions have the funny
rank one property. It would be interesting to see whether it is
possible to carry out Guenais's construction of ``good'' Morse
cocycles over a weakly mixing $(C,F)$-action of the group $G$ to
obtain a system which is not AT; as Morse cocycles over
$(C,F)$-actions yield systems which have simple spectrum on $L_1$
we have additionally to know that the spectral types on $L_0$ and
$L_1$ are mutually singular. Such a construction would be an
analog of Ageev's constructions \cite{Ag} for $\Z$-actions. It
would be even more interesting if we could carry out the above
over mixing $(C,F)$-actions of the group $G$ constructed by
Danilenko in \cite{Dan}.
\end{rem}

\begin{thank}\em
The  authors would like to  thank  J.-P. Thouvenot and F. Parreau
for fruitful discussions on the subject.
\end{thank}

\bibliographystyle{amsplain}
\setlength{\parsep}{0cm} \small
 %\bibliography{xbib}

\end{document}